\pdfoutput=1
\RequirePackage{ifpdf}
\ifpdf 
\documentclass[pdftex]{sigma}
\else
\documentclass{sigma}
\fi

\numberwithin{equation}{section}
\numberwithin{theorem}{section}
\numberwithin{proposition}{section}
\numberwithin{lemma}{section}
\numberwithin{corollary}{section}
\numberwithin{definition}{section}
\numberwithin{example}{section}
\numberwithin{remark}{section}
\numberwithin{note}{section}

\DeclareFontFamily{U}{mathb}{\hyphenchar\font45}
\DeclareFontShape{U}{mathb}{m}{n}{
      <5> <6> <7> <8> <9> <10> gen * mathb
      <10.95> mathb10 <12> <14.4> <17.28> <20.74> <24.88> mathb12
      }{}
\DeclareSymbolFont{mathb}{U}{mathb}{m}{n}
\DeclareFontSubstitution{U}{mathb}{m}{n}
\let\dot\relax
\DeclareMathAccent{\dot}{0}{mathb}{"39}
\let\ddot\relax
\DeclareMathAccent{\ddot}{0}{mathb}{"3A}
\let\dddot\relax
\DeclareMathAccent{\dddot}{0}{mathb}{"3B}
\let\ddddot\relax
\DeclareMathAccent{\ddddot}{0}{mathb}{"3C}

\begin{document}

\newcommand{\pd}[2]{\frac{\partial #1}{\partial #2}}

\newcommand{\arXivNumber}{1407.7919}

\allowdisplaybreaks

\renewcommand{\PaperNumber}{102}

\FirstPageHeading

\ShortArticleName{Particle Motion in Monopoles and Geodesics on Cones}

\ArticleName{Particle Motion in Monopoles and Geodesics on Cones}

\Author{Maxence MAYRAND}

\AuthorNameForHeading{M.~Mayrand}

\Address{Department of Mathematics and Statistics, McGill University,\\
805 Sherbrooke Street West, Montreal, Quebec, Canada, H3A 0B9}
\Email{\href{mailto:maxence.mayrand@mail.mcgill.ca}{maxence.mayrand@mail.mcgill.ca}}

\ArticleDates{Received July 31, 2014, in f\/inal form November 01, 2014; Published online November 04, 2014}

\Abstract{The equations of motion of a~charged particle in the f\/ield of Yang's $\mathrm{SU}(2)$ monopole in
5-dimensional Euclidean space are derived by applying the Kaluza--Klein formalism to the principal bundle
$\mathbb{R}^8\setminus\{0\}\to\mathbb{R}^5\setminus\{0\}$ obtained by radially extending the Hopf f\/ibration $S^7\to
S^4$, and solved by elementary methods.
The main result is that for every particle trajectory $\mathbf{r}:I\to\mathbb{R}^5\setminus\{0\}$, there is
a~4-dimensional cone with vertex at the origin on which $\mathbf{r}$ is a~geodesic.
We give an explicit expression of the cone for any initial conditions.}

\Keywords{particle motion; monopoles; geodesics; cones}

\Classification{70H06; 34A26; 53B50}

\section{Introduction}
The problem of the classical motion of an electrically charged particle in the f\/ield of Dirac's magnetic monopole is
a~system of three second-order non-linear dif\/ferential equations, written concisely as
\begin{gather}
\label{introEq}
\ddot{\mathbf{r}}=\lambda\frac{\mathbf{r}\times\dot{\mathbf{r}}}{|\mathbf{r}|^3},
\end{gather}
for $\mathbf{r}\in\dot{\mathbb{R}}^3:=\mathbb{R}^3\setminus\{0\}$ and a~constant $\lambda\in\mathbb{R}$.
We f\/ind it remarkable that, although Dirac's original paper~\cite{Dirac} about his monopole only appeared in 1931, Henri
Poincar\'e investigated the exact same system of equations in a~1896 paper~\cite{Poincare}.
His analysis was a~successful attempt to explain an experiment of the physicist Kristian Birkeland, which consisted of
approaching one pole of a~strong magnet near cathode rays, the other pole being far enough to be considered negligible.
We thus call this one-body dynamical system the ``Poincar\'e problem in $\dot{\mathbb{R}}^3$''.

In this paper, we are interested in the generalization of this problem to $\mathrm{SU}(2)$ gauge theory.
Recall that Dirac's monopole in $\mathbb{R}^3$ is obtained by radially extending the Hopf f\/ibration $S^3\to S^2$ to
a~principal $\mathrm{U}(1)$-bundle over $\dot{\mathbb{R}}^3$~\cite{Minami79,Ryder,Trautman}.
The same procedure using the next Hopf f\/ibration $S^7\to S^4$ gives rise to non-Abelian analogue of the monopole in
Euclidean space $\mathbb{R}^5$, known in the literature as Yang's monopole~\cite{Minami80,Yang}.
It is a~non-trivial $\mathrm{SO}(5)$-symmetrical solution to the Yang--Mills equations in $\mathrm{SU}(2)$ gauge theory.
Our main concern, which we call the ``Poincar\'e problem in $\dot{\mathbb{R}}^5$'', is for the classical motion of
a~charged particle in the presence of this monopole.
The equations of motion are derived in Section~\ref{Section3} using a~Kaluza--Klein formalism.
In this context, the charge~-- which generalizes~$\lambda$ in~\eqref{introEq}~-- is a~vector $\mathbf{e}$ rotating in
$\mathbb{R}^3$.

The f\/irst system,~\eqref{introEq}, has been thoroughly studied in the
literature~\cite{Feher,Fierz,Goddard,Haas,Horvathy,Jackiw,Lapidus,Moreira,Poincare,Ritter,Sivardiere}.
The main result~-- as shown f\/irst by Poincar\'e~-- is that for every solution $\mathbf{r}$, there is a~cone with vertex
at the origin on which $\mathbf{r}$ is a~geodesic (it follows from~\eqref{introEq} that $|\dot{\mathbf{r}}|$ is
constant).
Moreover, Poincar\'e provided an explicit expression for the cone's direction and the angle at its vertex (which vary
depending on the initial conditions and the charge~$\lambda$).
Since geodesics on cones are well understood, we get a~complete description of the space of solutions.

The main result of this paper is that this correspondence with geodesics on cones also holds for Yang's monopole, with
suitable modif\/ications.
Given any solution $\mathbf{r}:I\to\dot{\mathbb{R}}^5$ of the equations of motion, there is a~4-dimensional cone with
vertex at the origin of $\mathbb{R}^5$ on which $\mathbf{r}$ is a~geodesic.

Our proof proceeds in two main steps.
The f\/irst is the derivation of an explicit expression (given in Theorem~\ref{TheoremYangMonopole}) for the direction
$\mathbf{L}\in\dot{\mathbb{R}}^5$ of the 4-dimensional cone on which the particle is a~geodesic.
The second is a~general result (Theorem~\ref{TheoremConeGeo}) about geodesics on higher dimensional cones that we prove
here.
This theorem states that for all $n\geq 2$, a~geodesic on an~$n$-dimensional cone~$C$ is also a~geodesic on
a~$2$-dimensional cone embedded in~$C$ with the same angle at the vertex, and conversely.

Moreover, this last result shows that particles in Yang's monopole follow geodesics on 2-dimensional cones, and hence
all solutions can be obtained explicitly, as was the case for Dirac's monopole.

There is a~closely related problem called the ``MICZ-Kepler system''~\cite{MIC,Zwanzige}, which comes from generalizing
the Kepler problem (for the motion of a~particle under a~central inverse-squared attractive force in
$\dot{\mathbb{R}}^3$) by adding a~Lorentz force due to Dirac's monopole at the origin.
It has also been generalized in $\dot{\mathbb{R}}^5$ using Yang's monopole~\cite{Iwai}, and in all Euclidean spaces
$\dot{\mathbb{R}}^{n}$ by a~construction due to Meng~\cite{Meng2007-1,Meng2007-2,Meng2013}.
It was shown~\cite{BaiMengWang} that for all odd dimensions, the solutions to these systems are all conics.
Moreover, Montgomery showed~\cite{Montgomery2013} that in any dimension, this system is equivalent to the classical
Kepler problem on a~cone (with no magnetic charge).
It is thus natural to expect that the magnetic monopole alone would yield straight lines on cones (geodesics).
Our paper shows that this is the case, at least for Dirac's and Yang's monopole.

The paper is organized as follows.
In Section~\ref{Section2}, we recall the classical treatment of the Poincar\'e problem in $\dot{\mathbb{R}}^3$.

\looseness=-1
In Section~\ref{Section3} we brief\/ly review the Kaluza--Klein formalism for the motion of a~charged particle
in a~Yang--Mills f\/ield~\cite{Cho,Harnad,Kerner,Orzalesi}.
For a~principal $G$-bundle $P\!\to\!M$ with connection, the Kaluza--Klein approach is to construct
a~particular $G$-invariant metric on~$P$ from a~metric on~$M$ and an Ad-invariant metric on $\mathfrak{g}$.
Then, projection on~$M$ of the geodesics on~$P$ def\/ines motion of charged particles in~$M$.
The analogue of the charge is a~vector rotating in $\mathfrak{g}$.
The goal of this section is to provide coordinate expressions for the equations of motion.
We note (see Montgomery~\cite{Montgomery}) that this formulation is equivalent to the ones used~by
Sternberg~\cite{Sternberg}, Weinstein~\cite{Weinstein}, and Wong~\cite{Wong}.

In Section~\ref{Section4} we describe the extended Hopf bundles endowed with connections that give the Poincar\'e
problem in $\dot{\mathbb{R}}^3$ and $\dot{\mathbb{R}}^5$.
They are obtained by radially extending the Hopf f\/ibrations $S^{2n-1}\to S^n$ for $n=2,4$ to f\/ibrations
$\dot{\mathbb{R}}^{2n}\to\dot{\mathbb{R}}^{n+1}$, and taking the connections corresponding to a~horizontal subspace that
is orthogonal to the vertical subspace in Euclidean space $\dot{\mathbb{R}}^{2n}$.
As a~f\/irst example we apply the Kaluza--Klein formalism to $\dot{\mathbb{R}}^4\to\dot{\mathbb{R}}^3$ and show that we
recover the equations of motion~\eqref{introEq}.

In Section~\ref{Section6} we derive the equations of motion of the Poincar\'e problem in $\dot{\mathbb{R}}^5$.
That is, the one-body dynamical system for the motion of a~charged particle in the f\/ield of Yang's monopole.

Section~\ref{Section7} is devoted to the study of geodesics on higher dimensional cones.
This section is independent from the rest of the paper, but its conclusions will be crucial to the solution of the
Poincar\'e problem in $\dot{\mathbb{R}}^5$.

Finally, in Section~\ref{Section8} we show that a~charged particle in Yang's monopole must follow a~geodesic on
a~4-dimensional cone centred at the origin of $\mathbb{R}^5$.
We give an explicit expression for the cone, and thus obtain a~complete description of the space of solutions.

As a~side remark, we note that there is a~converse to the result of this paper.
We prove here one implication, namely, that if $\mathbf{r}$ a~solution to the Poincar\'e problem (in
$\dot{\mathbb{R}}^3$ or $\dot{\mathbb{R}}^5$) then $\mathbf{r}$ is a~geodesic on a~cone with vertex at the origin.
But we also have that for any cone centred at the origin (of $\mathbb{R}^3$ or $\mathbb{R}^5$) and any geodesic
$\mathbf{r}$ on it, there is a~unique charge ($\lambda$ or $\mathbf{e}$) for which $\mathbf{r}$ is a~solution to the
Poincar\'e problem (in $\dot{\mathbb{R}}^3$ or $\dot{\mathbb{R}}^5$).
For brevity we will not discuss this, but it can be proved with the theory presented in this paper.

\section{Particle motion in Dirac's monopole}\label{Section2}

Let us recall Poincar\'e's work~\cite{Poincare} on the motion of a~charged particle in the f\/ield of a~single magnetic pole.
Taking the pole to be centred at the origin, we f\/ind an electromagnetic f\/ield of the form
\begin{gather*}
\mathbf{E}=0,
\qquad
\mathbf{B}=g\frac{\mathbf{r}}{r^3},
\end{gather*}
for some constant $g\in\mathbb{R}$, $\mathbf{r}\in\dot{\mathbb{R}}^3:=\mathbb{R}^3\setminus\{0\}$ and $r=|\mathbf{r}|$.
Assuming the particle is subject to the Lorentz force $\mathbf{F}=q(\mathbf{E}+\dot{\mathbf{r}}\times \mathbf{B})$, we
get the equation of motion
\begin{gather}
\label{PoincareEquations}
\ddot{\mathbf{r}}=\lambda\frac{\mathbf{r}\times\dot{\mathbf{r}}}{r^3},
\end{gather}
for some constant $\lambda\in\mathbb{R}$ depending on the strength~$g$ of the magnet and the mass~$m$ and charge~$q$ of
the particle.
This is the system of ordinary dif\/ferential equations that Poincar\'e analysed, and is also the one describing motion of
a~charged particle in the f\/ield of Dirac's monopole.
Now, as Poincar\'e noticed, dif\/ferentiation shows that the vector
\begin{gather*}
\mathbf{L}:= \mathbf{r}\times\dot{\mathbf{r}}+\lambda\frac{\mathbf{r}}{r}
\end{gather*}
is constant.
Taking the norm, we see that $\mathbf{L}=0$ if and only if $\lambda=0$ and $\dot{\mathbf{r}}$ is everywhere parallel to
$\mathbf{r}$.
This corresponds to motion at constant speed on a~straight line that passes through the origin.
Since those curves will come out often here and in subsequent sections, we give them the following name (a term borrowed
from~\cite{BaiMengWang}).

\begin{definition}
A~\emph{colliding curve} is a~curve $\mathbf{r}:I\to\dot{\mathbb{R}}^n$ such that $\dot{\mathbf{r}}$ is everywhere
parallel to~$\mathbf{r}$.
\end{definition}

Now, suppose $\mathbf{r}$ is non-colliding.
Then, the cosine of the angle between $\mathbf{r}$ and $\mathbf{L}$ is
\begin{gather}
\label{angle}
\cos\psi=\frac{\mathbf{r}\cdot\mathbf{L}}{|\mathbf{r}||\mathbf{L}|}=\frac{\lambda}{|\mathbf{L}|},
\end{gather}
which is constant.
Hence, the particle moves on a~cone directed along $\mathbf{L}$.
Furthermore,~\eqref{PoincareEquations} shows that the acceleration is always normal to the surface of the cone, and so
the particle follows a~geodesic of that cone.

\begin{figure}[t]
\centering \includegraphics[scale=0.2]{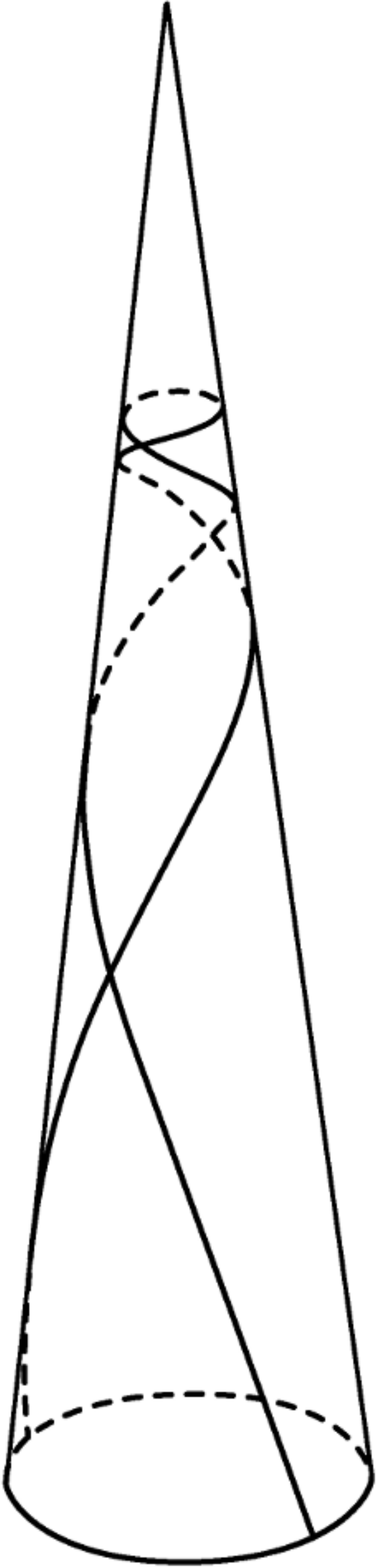}
\caption{A geodesic on a~cone.}\label{Fig1}
\end{figure}

With this information in hand, the problem reduces to the geodesic equations on a~cone in~$\mathbb{R}^3$~-- a standard
problem.
Note that the system~\eqref{PoincareEquations} is invariant under rotation, so we may assume the cone is directed along
the positive~$z$-axis.
Taking $t=0$ to be the point of closest approach to the origin, we f\/ind
\begin{gather*}
\mathbf{r}(t)=\sqrt{r_0^2+v_0^2t^2}\left(\sin\psi\cos\left(\frac{\arctan(v_0t/r_0)}{\sin\psi}\right),
\sin\psi\sin\left(\frac{\arctan(v_0t/r_0)}{\sin\psi}\right),\cos\psi\right),
\end{gather*}
where $r_0$, $v_0$ are the initial radius and velocity and~$\psi$ is half the angle at the vertex of the cone.
Moreover, equation~\eqref{angle} gives an explicit expression for the angle, namely, $\psi=\arctan(r_0 v_0/\lambda)$.
Fig.~\ref{Fig1} shows a~geodesic on a~cone.

In this section, the equations of motion were derived classically by considering the ``Coulomb-like'' magnetic f\/ield
$\mathbf{B}=g\mathbf{r}/|\mathbf{r}|^3$ and the Lorentz force.
But Dirac's monopole is also naturally described in terms of a~connection on the radial extension of the Hopf bundle
$S^1\to S^3\to S^2$.
In Section~\ref{Section4} we will show that this approach together with the Kaluza--Klein formalism give rise to the
exact same equations of motion.
Motion in Yang's monopole in $\dot{\mathbb{R}}^5$ will be obtained this way but by using the next Hopf f\/ibration $S^3\to
S^7\to S^4$.

\section{The Kaluza--Klein formalism}
\label{Section3}

Since Dirac's and Yang's monopole are more generally Yang--Mills f\/ields, we need a~way of obtaining the equations of
motion of a~particle in a~general Yang--Mills f\/ield.
There are several equivalent ways~\cite{Montgomery} of doing this, including the formulations of
Sternberg~\cite{Sternberg}, Weinstein~\cite{Weinstein}, Wong~\cite{Wong}, and Kerner~\cite{Kerner}.
In this paper we use the latter approach, which is known in the literature as the ``Kaluza--Klein formalism''.
The goal of this section is to brief\/ly review this formalism and to give coordinate expressions for the equations of
motion.
We closely follow the presentation of~\cite{Harnad}.
See also~\cite{Cho,Kerner,Orzalesi}.

Let $\pi:P\rightarrow M$ be a~principal bundle with structure group~$G$ acting on~$P$ to the right.
We then have local sections $\sigma_i:U_i\to\pi^{-1}(U_i)$ such that the local trivialisations $\phi_i(p)=(\pi(p),a)\in
U_i\times G$ correspond to the right action $p=\sigma(\pi(p))\cdot a$.
Let~$\theta$ be a~connection one-form on~$P$.
For an Ad-invariant metric $\langle~,\,\rangle$ on the Lie algebra $\mathfrak{g}$ of~$G$ and a~metric~$g$ on~$M$, def\/ine
the following $G$-invariant metric on~$P$,
\begin{gather*}
\gamma(X,Y)|_p=g(\pi_\ast(X),\pi_\ast(Y))|_{\pi(p)}+\langle\theta(X)|_p,\theta(Y)|_p\rangle.
\end{gather*}
Then,~$P$ with this metric~$\gamma$ is a~Riemmanian manifold whose geodesics projected to~$M$ def\/ine the motion of
charged particles in~$M$, where the charge is a~vector of constant magnitude in $\mathfrak{g}$.

We now set up the equations of motion in terms of a~local coordinate system $\{x_i\}$ on an open neighbourhood
$U\subseteq M$ and a~basis $\{T_k\}$ for $\mathfrak{g}$.
Let $T_pP=V_p\oplus H_p$ be the decomposition of the tangent space into a~vertical and horizontal subspace.
The action of~$G$ induces a~canonical isomorphism between $V_p$ and $\mathfrak{g}$, which gives fundamental vector
f\/ields $\{L_k\}$ on~$P$ corresponding to $\{T_k\}$.
Then, for a~curve $t\mapsto p(t)$ in~$P$, the tangent vector at~$t$ is $\dot{p}(t)=v(t)+h(t)$, for some $v(t)\in
V_{p(t)}$ and $h(t)\in H_{p(t)}$, and we may expand $v=v^kL_k$.
The geodesic equations for the curve~$p$ with respect to the metric~$\gamma$ then become
\begin{gather}
\ddot{x}^\mu+\Gamma^\mu_{\lambda\rho}\dot{x}^\lambda\dot{x}^\rho = \langle v,\tilde{F}_{\lambda\nu}\dot{x}^\nu\rangle g^{\lambda\mu},
\label{geodesic1}
\\
\dot{v}^k =0,
\label{geodesic2}
\end{gather}
where $\Gamma^\mu_{\lambda\rho}$ are the Christof\/fel symbols of the metric~$g$ on~$M$, and
\begin{gather*}
\tilde{F}=d\theta+\frac{1}{2}[\theta,\theta]= \frac{1}{2}\tilde{F}_{\mu\nu}\mathrm{d} x^\mu\wedge\mathrm{d} x^\nu
\end{gather*}
is the curvature two-form.
Let $\sigma:U\to\pi^{-1}(U)$ be the canonical local section and def\/ine the $\mathfrak{g}$-valued one-form
\begin{gather*}
A=A_\mu\mathrm{d} x^\mu=A^k_\mu T_k\mathrm{d} x^\mu=\sigma^\ast(\theta),
\end{gather*}
called the local gauge potential over~$U$.
Let $a\in G$ be the local trivialization $p\sim(\pi(p),a)\in U\times G$, and let
\begin{gather*}
e=e^kT_k=ava^{-1}.
\end{gather*}
Finally, set $b_{ij}:=\langle T_i,T_j\rangle$, let $C^k_{ij}$ be the structure constants of $\mathfrak{g}$, and let
\begin{gather*}
F_{\mu\nu}=F^k_{\mu\nu}T_k=\partial_\mu A_\nu-\partial_\nu A_\mu+[A_\mu,A_\nu]
\end{gather*}
be the curvature two-form in the~$\sigma$-gauge.
So $F_{\mu\nu}$ and $\tilde{F}_{\mu\nu}$ are related by $\tilde{F}_{\mu\nu}=a^{-1} F_{\mu\nu}a$.
Then, the geodesic equations~\eqref{geodesic1} and~\eqref{geodesic2} are
\begin{gather}
\label{generalequationsofmotion1}
\ddot{x}^\mu+\Gamma^\mu_{\lambda\rho}\dot{x}^\lambda\dot{x}^\rho = b_{ij}e^iF^j_{\lambda\nu}\dot{x}^\nu g^{\lambda\mu},
\\
\label{generalequationsofmotion2}
\dot{e}^k+C^k_{ij}A^i_\mu e^j \dot{x}^\mu = 0.
\end{gather}
This system of ordinary dif\/ferential equations def\/ines the motion of a~charged particle in~$M$.
The term on the right-hand side of~\eqref{generalequationsofmotion1} is the generalization of the Lorentz force.
The vector $e=e^kT_k\in\mathfrak{g}$ is the analogue of the charge divided by the mass of the particle.
Note that~$e$ has magnitude $\langle e,e\rangle^{1/2}=\langle v,v\rangle^{1/2}$, which is constant by~\eqref{geodesic2}.
However, unless~$G$ is Abelian, $e$~itself is not in general constant.

\section{The extended Hopf bundles}
\label{Section4}

The purpose of this section is to def\/ine the principal bundles endowed with connections that describe Dirac's and Yang's
monopole.
The bundles are obtained by radially extending the Hopf bundles $S^{n-1}\to S^{2n-1}\to S^n$, for $n=2,4$.
We will f\/irst give a~more abstract def\/inition by means of the canonical projection of certain quotient spaces of the
vector spaces $\dot{\mathbb{C}}^2$ and $\dot{\mathbb{H}}^2$.
It will then lead to the desired bundles $\dot{\mathbb{R}}^{2n}\to\dot{\mathbb{R}}^{n+1}$ by dif\/feomorphisms.
Similar constructions can be found in~\cite{Trautman} and~\cite{Duval}.

Let~$K$ be $\mathbb{C}$ or $\mathbb{H}$, and let $\dot{K}^2:=K^2\setminus\{(0,0)\}$.
Let $\sim$ be the equivalence relation on $\dot{K}^2$ def\/ined by $(z_1,z_2)\sim(w_1,w_2)$ if there is a~unit norm
$\lambda\in K$ such that $(z_1,z_2)=(w_1\lambda,w_2\lambda)$.
The quotient of $\dot{K}^2$ by this relation gives an $(n+1)$-dimensional dif\/ferentiable manifold~$M$, and we def\/ine the
extended Hopf map by the canonical projection
\begin{gather}
\label{AbstractExtendedHopfFibration}
\pi: \ \dot{K}^2\to M.
\end{gather}
We get the structure of a~principal bundle as follows.
Let $S^1_K$ be the set of unit norm elements in~$K$.
We have $S^1_\mathbb{C}=\mathrm{U}(1)$ and $S^1_\mathbb{H}=\mathrm{SU}(2)$, so $S^1_K$ is a~Lie group.
It acts freely on $\dot{K}^2$ by $(z_1,z_2)\cdot\lambda=(z_1\lambda,z_2\lambda)$, and~$M$ is the quotient space of this
action.
Moreover,~$M$ is covered by the two open neighbourhoods
\begin{gather*}
U_i:=\{[(z_1,z_2)]\in M:z_i\neq 0\},
\qquad
i=1,2,
\end{gather*}
over which we have the local trivializations
\begin{gather}
\label{LocalTriv}
\pi^{-1}(U_i)\to M\times S^1_K,
\qquad
(z_1,z_2)\mapsto ([(z_1,z_1)],z_i/|z_i|),
\qquad
i=1,2.
\end{gather}
Thus, the extended Hopf map~\eqref{AbstractExtendedHopfFibration} is a~principal $S^1_K$-bundle.
Now, the principal bundle $S^1_K\to\dot{\mathbb{R}}^{2n}\to\dot{\mathbb{R}}^{n+1}$ is obtained
from~\eqref{AbstractExtendedHopfFibration} by the identif\/ication of~$K$ with $\mathbb{R}^n$ using the basis $\{1,i\}$
for~$\mathbb{C}$ and $\{1,i,j,k\}$ for $\mathbb{H}$, and by the dif\/feomorphism
\begin{gather*}
f: \ M\to\dot{\mathbb{R}}^{n+1},
\qquad
[(z_1,z_2)]\mapsto\left(\frac{2z_1z_2^*}{\sqrt{|z_1|^2+|z_2|^2}},\frac{|z_1|^2-|z_2|^2}{\sqrt{|z_1|^2+|z_2|^2}}\right),
\end{gather*}
where $^*$ denotes conjugation in~$K$.
This construction gives the two principal bundles
\begin{gather*}
\mathrm{U}(1)\to \dot{\mathbb{R}}^4\to\dot{\mathbb{R}}^3
\qquad
\text{and}
\qquad
\mathrm{SU}(2)\to\dot{\mathbb{R}}^8\to\dot{\mathbb{R}}^5.
\end{gather*}
Motion will take place in Euclidean spaces $\dot{\mathbb{R}}^3$ and $\dot{\mathbb{R}}^5$.

Let us introduce the following set of coordinates on $M\cong\dot{\mathbb{R}}^{n+1}$.
\begin{gather}
\phi_1: \ U_1\subseteq M\to\mathbb{R}^n\times\mathbb{R}^+,
\qquad
[(z_1,z_2)]\mapsto\big(z_2z_1^{-1},\sqrt{|z_1|^2+|z_2|^2}\big),
\\
\phi_2: \ U_2\subseteq M\to\mathbb{R}^n\times\mathbb{R}^+,
\qquad
[(z_1,z_2)]\mapsto \big(z_1z_2^{-1},\sqrt{|z_1|^2+|z_2|^2}\big).
\label{Coordinates}
\end{gather}
These coordinates are denoted $(\mathbf{u},r)=(u_1,\ldots,u_n,r)$ and are related to the Cartesian coordinates
$(x_1,\ldots,x_{n+1})$ of $\dot{\mathbb{R}}^{n+1}$~by
\begin{gather}
(x_1,\ldots,x_{n+1}) =
\big(f\circ\phi_1^{-1}\big)(\mathbf{u},r)=\left(\frac{2r\mathbf{u}^*}{1+|\mathbf{u}|^2},r\frac{1-|\mathbf{u}|^2}{1+|\mathbf{u}|^2}\right),
\nonumber
\\
(x_1,\ldots,x_{n+1}) =
\big(f\circ\phi_2^{-1}\big)(\mathbf{u},r)=\left(\frac{2r\mathbf{u}}{|\mathbf{u}|^2+1},r\frac{|\mathbf{u}|^2-1}{|\mathbf{u}|^2+1}\right).
\label{StereoNorth}
\end{gather}
An observation that will be crucial later is that $(\mathbf{u},r)$ are precisely the stereographic projection
coordinates from the south and north poles respectively.
That is, for $\mathbf{r}\in\dot{\mathbb{R}}^{n+1}$, f\/irst project on the unit sphere by $\mathbf{r}\mapsto
\mathbf{r}/|\mathbf{r}|$.
Then, the stereographic projection of $\mathbf{r}/|\mathbf{r}|$ gives $\mathbf{u}=(u_1,\ldots,u_n)$, and the remaining
coordinate~$r$ is the magnitude of $\mathbf{r}$.

To apply the Kaluza--Klein formalism, we further need a~connection on the bundle, a~metric on
$M\cong\dot{\mathbb{R}}^{n+1}$ and an Ad-invariant metric on $\mathfrak{g}$.
The connection that gives the Poincar\'e problem is obtained by choosing a~horizontal subspace that is orthogonal to the
vertical subspace in Euclidean space $\dot{\mathbb{R}}^{2n}$.
The metric on~$M$ is the one corresponding to the Euclidean metric on $\dot{\mathbb{R}}^{n+1}$.
Since we observed that the coordinates $(u_1,\ldots,u_n,r)$ on~$M$ are the stereographic projection coordinates of
$\dot{\mathbb{R}}^{n+1}$, we know as a~standard result that the metric on~$M$ is
\begin{gather*}
g=\frac{4r^2\sum\limits_{i=1}^n\mathrm{d} u_i\otimes\mathrm{d} u_i}{(1+u_1^2+\dots+u_n^2)^2}+\mathrm{d}
r\otimes\mathrm{d} r.
\end{gather*}
For the Ad-invariant metric on $\mathfrak{g}$ we take $\langle T_i,T_j\rangle:=\delta_{ij}$, where $\{T_i\}$ is a~basis
for $\mathfrak{g}$.
When $G=\mathrm{U}(1)$, this basis is the imaginary number $T_1=i$, and when $G=\mathrm{SU}(2)$, it is
$\{T_1=i$, $T_2=j,T_3=k\}$.
It is straightforward to verify that $\langle T_i,T_j\rangle:=\delta_{ij}$ is Ad-invariant.

As an example, we apply the Kaluza--Klein formalism to the principal bundle
$\mathrm{U}(1)\to\dot{\mathbb{R}}^4\to\dot{\mathbb{R}}^3$ and show that we recover the equations of motion obtained in
Section~\ref{Section2}, i.e.~those describing the classical motion of a~charged particle in the f\/ield of Dirac's monopole.

The vertical subspace $V_p\subseteq T_p\dot{\mathbb{C}}^2$ for $p=(z_1,z_2)\in \dot{\mathbb{C}}^2$ is spanned~by
$\frac{d}{dt}\big|_{t=0}(z_1,z_2)\cdot\exp(ti)=(z_1i,z_2i)$.
In the Cartesian coordinates $(x_1,x_2,x_3,x_4)=(z_1,z_2)$ of $\dot{\mathbb{R}}^4$, we have
\begin{gather*}
V_p=\mathrm{span}\left\{-x_2\pd{}{x_1}+x_1\pd{}{x_2}-x_4\pd{}{x_3}+x_3\pd{}{x_4}\right\}.
\end{gather*}
We take the horizontal subspace $H_p$ to be the orthogonal complement of $V_p$.
The corresponding connection one-form is then
\begin{gather*}
\theta=i\frac{-x_2\mathrm{d} x_1+x_1\mathrm{d} x_2-x_4\mathrm{d} x_3+x_3\mathrm{d} x_4}{x_1^2+\dots+x_4^3}.
\end{gather*}
This gives the local gauge potential over $U_1$
\begin{gather*}
A=\sigma_1^*(\theta)=i\frac{ydx-xdy}{2r(z+r)},
\end{gather*}
where $(x,y,z)$ are the Cartesian coordinates in $\dot{\mathbb{R}}^3$ and $r=\sqrt{x^2+y^2+z^2}$.
The curvature then reads
\begin{gather*}
F=i\frac{-z
dx\wedge dy+y
dx\wedge dz-x
dy\wedge dz}{2r^3}.
\end{gather*}
Inserting in the equations of motion~\eqref{generalequationsofmotion1} and~\eqref{generalequationsofmotion2}, we get
\begin{gather*}
\ddot{x}=\frac{e}{2}\cdot\frac{y\dot{z}-z\dot{y}}{r^3},
\qquad
\ddot{y}=\frac{e}{2}\cdot\frac{z\dot{x}-x\dot{z}}{r^3},
\qquad
\ddot{z}=\frac{e}{2}\cdot\frac{x\dot{y}-y\dot{x}}{r^3},
\qquad
\dot{e}=0,
\end{gather*}
which are precisely Poincar\'e's original equations~\eqref{PoincareEquations} with $\lambda=e/2$.

\section{The equations of motion of a~particle in Yang's monopole}
\label{Section6}

In this section, we obtain the equations of motion of a~charged particle in the presence of Yang's monopole in
$\dot{\mathbb{R}}^5$ by applying the Kaluza--Klein formalism to the extended Hopf bundle
$\mathrm{SU}(2)\to\dot{\mathbb{R}}^8\to\dot{\mathbb{R}}^5$ constructed in Section~\ref{Section4}.

To compute the vertical subspace, we use the right action of $\mathrm{SU}(2)$ to pullback the basis $\{i,j,k\}$ of
$\mathfrak{su}(2)$ to a~basis for $V_p$.
The basis vectors are $L_i|_p:=\frac{d}{dt}\big|_{t=0} p\cdot\exp(tT_i)$, and in the Cartesian coordinates
$(x_1,\ldots,x_8)=(z_1,z_2)$ of $\dot{\mathbb{R}}^8\cong\dot{\mathbb{H}}^2$, we have
\begin{gather*}
  L_1  =  (-x_2,  x_1,  x_4,  -x_3,  -x_6,  x_5,  x_8,  -x_7),
\\
L_2  =  (-x_3,  -x_4,  x_1,  x_2,  -x_7,  -x_8,  x_5,  x_6),
\\
L_3  =  (-x_4,  x_3,  -x_2,  x_1,  -x_8,  x_7,  -x_6,  x_5) .
\end{gather*}
The connection one-form corresponding to a~horizontal subspace that is orthogonal to $V_p$ is then
\begin{gather*}
\theta=\frac{1}{x_1^2+\dots+x_8^2}
\begin{pmatrix}
-x_2&x_1&x_4&-x_3&-x_6&x_5&x_8&-x_7
\\
-x_3&-x_4&x_1&x_2&-x_7&-x_8&x_5&x_6
\\
-x_4&x_3&-x_2&x_1&-x_8&x_7&-x_6&x_5
\end{pmatrix}
\begin{pmatrix}
\mathrm{d} x_1
\\
\vdots
\\
\mathrm{d} x_8
\end{pmatrix},
\end{gather*}
as expressed in the basis $\{i,j,k\}$ for $\mathfrak{su}(2)$.
We will now use the coordinate system $(\mathbf{u},r)=(u_1,\ldots,u_4,r)$ on $U_2\subseteq M\cong\dot{\mathbb{R}}^5$
def\/ined by~\eqref{Coordinates}.
Recall that these are the stereographic projection coordinates from the north pole.
Hence, we are working in $\dot{\mathbb{R}}^5$ minus the positive $x_5$-axis.
Let $\sigma_i$ for $i=1,2$ be the canonical local sections induced by the local trivializations~\eqref{LocalTriv}.
On $U_2$ we obtain the gauge potential
\begin{gather*}
A=\sigma_2^\ast(\theta)=\frac{\mathbf{u}^*\mathrm{d}\mathbf{u}-\mathrm{d}\mathbf{u}^*\mathbf{u}}{2(|\mathbf{u}|^2+1)},
\end{gather*}
and the corresponding curvature
\begin{gather*}
F=\frac{\mathrm{d}\mathbf{u}^*\wedge\mathrm{d}\mathbf{u}}{2(|\mathbf{u}|^2+1)^2}.
\end{gather*}
In matrix notation, we have $A=\mathbf{A} \mathrm{d}\mathbf{u}$, where
\begin{gather*}
\mathbf{A}:= \frac{1}{|\mathbf{u}|^2+1}
\begin{pmatrix}
-u_2 & u_1 & u_4 & -u_3
\\
-u_3 & -u_4 & u_1 & u_2
\\
-u_4 & u_3 & -u_2 & u_1
\end{pmatrix}
=:
\begin{pmatrix}
\mathbf{A}_1
\\
\mathbf{A}_2
\\
\mathbf{A}_3
\end{pmatrix}
.
\end{gather*}
Inserting in the equations of motions~\eqref{generalequationsofmotion1} and~\eqref{generalequationsofmotion2} of the
Kaluza--Klein formalism, we get the system of dif\/ferential equations
\begin{gather}
\label{SU2EOM1}
\ddot{\mathbf{u}}+\frac{2|\dot{\mathbf{u}}|^2\mathbf{u}-4(\mathbf{u}\cdot\dot{\mathbf{u}})\dot{\mathbf{u}}}{|\mathbf{u}|^2+1}
+\frac{2\dot{r}\dot{\mathbf{u}}}{r}
 = \frac{\mathbf{E}\dot{\mathbf{u}}}{2r^2},
\\
\label{SU2EOM2}
\ddot{r}-\frac{4r|\dot{\mathbf{u}}|^2}{(|\mathbf{u}|^2+1)^2} = 0,
\\
\label{SU2EOM3}
\dot{\mathbf{e}}+2\mathbf{B}\mathbf{e} = 0,
\end{gather}
where $\mathbf{e}=(e^1,e^2,e^3)$ and
\begin{gather*}
\mathbf{E}:=
\begin{pmatrix}
0&e^1&e^2&e^3
\\
-e^1&0&-e^3&e^2
\\
-e^2&e^3&0&-e^1
\\
-e^3&-e^2&e^1&0,
\end{pmatrix}
,
\qquad
\mathbf{B}:=
\begin{pmatrix}
0 & -B_3 & B_2
\\
B_3 & 0 & -B_1
\\
-B_2 & B_1 & 0
\end{pmatrix}
,
\qquad
B_i:=\mathbf{A}_i\cdot\dot{\mathbf{u}}.
\end{gather*}
These equations describe the motion of a~charged particle in the f\/ield of Yang's monopole at the origin of Euclidean
space $\mathbb{R}^5$.
The vector $\mathbf{e}$ is interpreted as the charge of the particle, and
$(\frac{1}{2r^2}\mathbf{E}\dot{\mathbf{u}},0)\in\mathbb{R}^4\times\mathbb{R}^+$ is the analogue of the Lorentz force.
Note that~\eqref{SU2EOM3} immediately gives $\mathbf{e}\cdot\dot{\mathbf{e}}=0$, and so $\mathbf{e}$ has constant
magnitude, as anticipated in the general formulation of Section~\ref{Section3}.

\section{Some facts about cones and their geodesics}
\label{Section7}

The solutions to the equations of motion~\eqref{SU2EOM1},~\eqref{SU2EOM2} and~\eqref{SU2EOM3} will be investigated in
Section~\ref{Section8}.
Some crucial results that we will need can be stated as general facts about higher dimensional cones and their
geodesics.
Hence, we put them in this separate section, which is completely independent from the rest of the paper.
The main goal is Theorem~\ref{TheoremConeGeo} and its two corollaries.

First of all, we need a~clear def\/inition of what we mean by a~$k$-dimensional cone in $\mathbb{R}^n$, for $k<n$.
In this paper all cones will have their vertex at the origin.
Before the general def\/inition, here is the most basic classical one.

\begin{definition}
\label{DefConeBasic}
The \emph{cone} of aperture $\psi\in(0,\pi/2]$ directed along $\mathbf{L}\in\dot{\mathbb{R}}^n$ is the set of all
points $\mathbf{r}\in\dot{\mathbb{R}}^n$ satisfying
\begin{gather}
\label{ConeDef}
\frac{\mathbf{r}\cdot\mathbf{L}}{|\mathbf{r}||\mathbf{L}|}=\cos\psi.
\end{gather}
We write a~``cone in $\mathbb{R}^n$'' or equivalently an ``$(n-1)$-dimensional cone in $\mathbb{R}^n$'' for any such
cone.
\end{definition}

\begin{figure}[t]
\centering \includegraphics[scale=0.2]{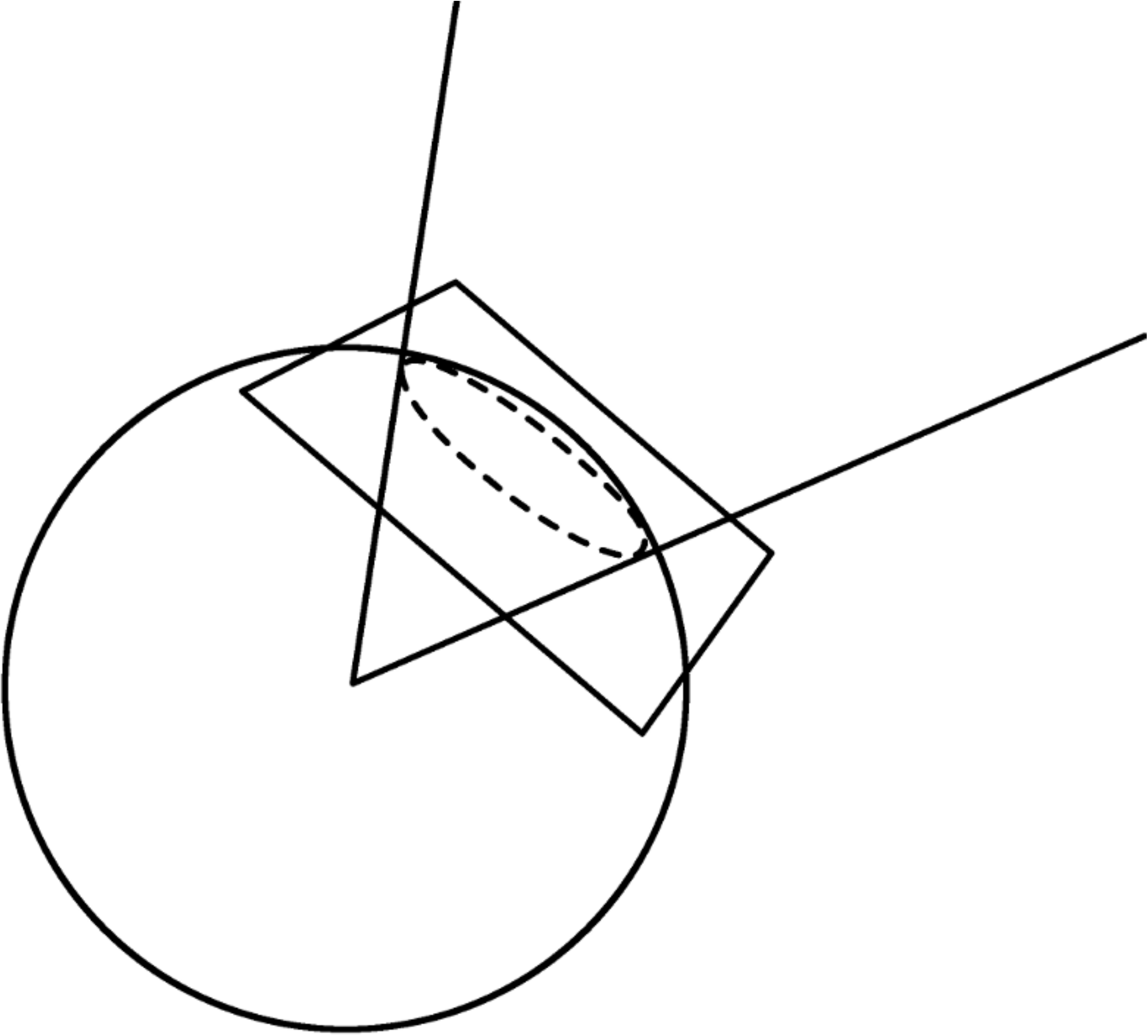}
\caption{A cone is def\/ined by radially extending the intersection of
an af\/f\/ine plane with the unit sphere.}\label{Fig2}
\end{figure}

To generalize this def\/inition to~$k$-dimensional cones in $\mathbb{R}^n$ for any $k<n$, we need the following
observation.
In the def\/inition of a~cone, equation~\eqref{ConeDef} can be rewritten $(\mathbf{r}/|\mathbf{r}|)\cdot\mathbf{L}=b$,
where $b=|\mathbf{L}|\cos\psi$ is a~constant.
Thus, the cone is the set of all points in $\mathbb{R}^n$ such that when projected on the unit sphere $S^{n-1}$ they lie
in the f\/ixed af\/f\/ine hyperplane $\{\mathbf{x}\in\mathbb{R}^n:\mathbf{x}\cdot\mathbf{L}=b\}$~-- see Fig.~\ref{Fig2}.
This motivates the following def\/inition.

\begin{definition}
\label{DefConeGen}
Let~$P$ be an af\/f\/ine~$k$-dimensional plane in $\mathbb{R}^n$ that intersects with $S^{n-1}$ in more than one point.
The \emph{cone generated by~$P$ in $\mathbb{R}^n$} is the set of all points $\mathbf{r}\in\dot{\mathbb{R}}^n$ such
that $\mathbf{r}/|\mathbf{r}|\in P$.
We write a~``$k$-dimensional cone in $\mathbb{R}^n$'' for any such set.
\end{definition}

To better understand this def\/inition, we f\/irst remark the following basic fact.

\begin{proposition}
\label{PlaneCapCircle}
Let~$P$ be a~$k$-dimensional affine plane in $\mathbb{R}^n$ intersecting the unit sphere $S^{n-1}$ in more than one
point.
Then, $P\cap S^{n-1}$ is a~$(k-1)$-sphere.
\end{proposition}

This proposition shows that a~$k$-dimensional cone in $\mathbb{R}^n$ is parametrized by a~point in a~$(k-1)$-sphere and
a~positive number $r>0$.
For example, a~$2$-dimensional cone in $\mathbb{R}^3$ intersects the unit sphere on a~circle, and so each point of the
cone is uniquely def\/ined by a~point on this circle and a~radius $r>0$.
Here is the proof of the proposition.

\begin{proof}
Let $P=\mathbf{a}+U$, where $\mathbf{a}\in\mathbb{R}^n$ and~$U$ is a~$k$-dimensional subspace of $\mathbb{R}^n$.
Without loss of generality, we may assume that $\mathbf{a}$ is orthogonal to~$U$.
Choose an orthonormal basis $\{\mathbf{v}_1,\ldots,\mathbf{v}_k\}$ for~$U$.
Then,
\begin{gather*}
P\cap S^{n-1}
=\big\{\mathbf{a}+x_1\mathbf{v}_1+\dots+x_k\mathbf{v}_k:x_i\in\mathbb{R},|\mathbf{a}+x_1\mathbf{v}_1+\dots+x_k\mathbf{v}_k|^2=1\big\}
\\
\phantom{P\cap S^{n-1}}
 = \big\{\mathbf{a}+x_1\mathbf{v}_1+\dots+x_k\mathbf{v}_k:x_i\in\mathbb{R},x_1^2+\dots+x_k^2=1-|\mathbf{a}|^2\big\},
\end{gather*}
which is a~$(k-1)$-sphere of radius $\sqrt{1-|\mathbf{a}|^2}$ centred at $\mathbf{a}$ in $\mathbb{R}^n$.
\end{proof}

The following proposition shows that a~$k$-dimensional cone in $\mathbb{R}^n$, as of Def\/inition~\ref{DefConeGen}, is in
a~sense exactly the same as the classical Def\/inition~\ref{DefConeBasic} of a~cone in $\mathbb{R}^{k+1}$.

\begin{proposition}
Let~$C$ be a~$k$-dimensional cone in $\mathbb{R}^n$ for $n>k$.
There is a~cone~$D$ in $\mathbb{R}^{k+1}$ directed along $(0,\ldots,0,1)\in\mathbb{R}^{k+1}$, and a~matrix
$\mathbf{R}\in\mathrm{SO}(n)$ such that
\begin{gather*}
\mathbf{R}(C)=\big\{(x_1,\ldots,x_{k+1},0,\ldots,0)\in\mathbb{R}^n:(x_1,\ldots,x_{k+1})\in D\big\}.
\end{gather*}
\end{proposition}

\begin{proof}
Let~$P$ be the~$k$-dimensional af\/f\/ine plane in $\mathbb{R}^n$ that generates~$C$.
It is straightforward to see that there exists $\mathbf{R}\in\mathrm{SO}(n)$ such that
\begin{gather*}
\mathbf{R}(P)=\{(x_1,\ldots,x_k,a,0,\ldots,0):x_i\in\mathbb{R}\},
\end{gather*}
for some $a\geq 0$.
Note that since~$P$ intersects $S^{n-1}$ in more than one point, we have $0\leq a<1$, and hence there is an angle
$\psi\in(0,\pi/2]$ such that $a=\cos\psi$.
Letting $\mathbf{L}=(0,\ldots,0,1)\in\mathbb{R}^{k+1}$, we get
\begin{gather*}
\mathbf{R}(C) = \big\{\mathbf{R}(\mathbf{r})\in\mathbb{R}^n:\mathbf{r}/|\mathbf{r}|\in P\big\} =
\big\{\mathbf{r}\in\mathbb{R}^n:\mathbf{r}/|\mathbf{r}|\in\mathbf{R}(P)\big\}
\\
\phantom{\mathbf{R}(C)}
 = \big\{\mathbf{r}\in\mathbb{R}^n:\mathbf{r}/|\mathbf{r}|=(x_1,\ldots,x_k,\cos\psi,0,\ldots,0)
\;
\text{for some}
\;
x_i\in\mathbb{R}\big\}
\\
\phantom{\mathbf{R}(C)}
=\left\{(\tilde{\mathbf{r}},0,\ldots,0)\in\mathbb{R}^{k+1}\times\mathbb{R}^{n-k-1}:
\frac{\tilde{\mathbf{r}}\cdot\mathbf{L}}{|\tilde{\mathbf{r}}||\mathbf{L}|}=\cos\psi\right\}
\\
\phantom{\mathbf{R}(C)}
=\big\{(x_1,\ldots,x_{k+1},0,\ldots,0)\in\mathbb{R}^n:(x_1,\ldots,x_{k+1})\in D\big\},
\end{gather*}
where~$D$ is the cone of aperture~$\psi$ directed along $\mathbf{L}$ in $\dot{\mathbb{R}}^{k+1}$.
\end{proof}

This proposition and its proof allow us to make the following def\/inition.

\begin{definition}
Let~$C$ be a~$k$-dimensional cone in $\mathbb{R}^n$, and let~$P$ be the af\/f\/ine~$k$-dimensional plane generating~$C$.
Write $P=\mathbf{a}+U$ for $\mathbf{a}\in U^\perp$.
The \emph{aperture} of~$C$ is the number $\psi\in(0,\pi/2]$ for which $\cos\psi=|\mathbf{a}|$.
\end{definition}

It is straightforward to verify that this is well-def\/ined ($\mathbf{a}$ is unique and $0\leq|\mathbf{a}|<1$) and that it
matches the classical Def\/inition~\ref{DefConeBasic}.
For a~2-dimensional cone~$C$ in $\mathbb{R}^3$, our def\/inition of the aperture is \emph{half} the angle at the vertex
of~$C$.

Consider a~1-dimensional cone of aperture~$\psi$ in $\mathbb{R}^3$ (two non-parallel rays coming from the origin).
It is intuitively clear that there is a~unique 2-dimensional cone of aperture~$\psi$ containing it.
Indeed, just rotate the two rays about the bisector, and it will give the desired cone.
This principle of ``unique embedding'' is indeed true, and generalizes as follows.

\begin{proposition}
\label{UniqueConeEmbeding}
Let~$D$ be a~$k$-dimensional cone of aperture~$\psi$ in $\mathbb{R}^{n}$, for any $k<n$.
There is a~unique $(n-1)$-dimensional cone of aperture~$\psi$ containing~$D$.
\end{proposition}

\begin{proof}
Let $Q=\mathbf{a}+U$ be the~$k$-dimensional af\/f\/ine plane generating~$D$, and assume $\mathbf{a}\in U^\perp$ so that
$\cos\psi=|\mathbf{a}|$.
Then, $U\subseteq(\mathrm{span}\{\mathbf{a}\})^\perp$, so the hyperplane
$P=\mathbf{a}+(\mathrm{span}\{\mathbf{a}\})^\perp$ generates an $(n-1)$-di\-men\-sional cone of aperture~$\psi$
containing~$D$.
This shows existence.
For uniqueness, suppose that $\tilde{P}=\mathbf{b}+V$, $\mathbf{b}\in V^\perp$, generates an $(n-1)$-di\-men\-sional cone of
aperture~$\psi$ containing~$D$.
Write $\mathbf{b}=\mathbf{a}+\mathbf{c}$ for $\mathbf{c}\in(\mathrm{span}\{\mathbf{a}\})^\perp$.
Then, $|\mathbf{a}|^2=\cos^2\psi=|\mathbf{b}|^2=|\mathbf{a}|^2+|\mathbf{c}|^2$, so $|\mathbf{c}|^2=0$ and hence
$\mathbf{b}=\mathbf{a}$.
Then, $\mathbf{a}\in V^\perp$ so we have $V\subseteq (\mathrm{span}\{\mathbf{a}\})^\perp$.
But $\dim V=n-1=\dim (\mathrm{span}\{\mathbf{a}\})^\perp$, so $V=(\mathrm{span}\{\mathbf{a}\})^\perp$, whence
$\tilde{P}=P$.
\end{proof}

We will now start to investigate geodesics on cones.
First, let us give a~clear def\/inition.

\begin{definition}
Let~$C$ be a~$k$-dimensional cone in $\mathbb{R}^n$ together with the metric inherited from the ambient Euclidean space
$\mathbb{R}^n$ and let $\frac{D}{dt}$ be the corresponding covariant derivative.
We call a~\emph{geodesic} on~$C$ a~dif\/ferentiable map $\mathbf{r}:I\to\dot{\mathbb{R}}^n$, from an open interval
$I\subseteq\mathbb{R}$ to $\dot{\mathbb{R}}^n$, such that $\mathbf{r}(I)\subseteq C$ and
$\frac{D}{dt}\dot{\mathbf{r}}(t)=0$ for all $t\in I$.
\end{definition}

Note that we do not assume that a~geodesic is parametrized by arclength, and hence can have any (constant) speed.

We will now need a~parametrization for an arbitrary~$n$-dimensional cone~$C$ in $\mathbb{R}^{n+1}$.
Since for any $\mathbf{R}\in\mathrm{SO}(n+1)$, a~curve $\mathbf{r}$ is a~geodesic on a~cone~$C$ if and only if
$\mathbf{R}(\mathbf{r})$ is a~geodesic on $\mathbf{R}(C)$, we may assume without loss of generality that~$C$ is directed
along $\mathbf{L}=(0,\ldots,0,1)\in\mathbb{R}^{n+1}$.
Let $\psi\in(0,\pi/2]$ be the aperture of~$C$.
Then,
\begin{gather*}
C = \Big\{(x_1,\ldots,x_{n+1})\in\mathbb{R}^{n+1}:x_{n+1}=\cos\psi\sqrt{x_1^2+\dots+x_{n+1}^2}\Big\}
\\
\phantom{C}
 = \big\{(x_1,\ldots,x_{n+1})\in\mathbb{R}^{n+1}:x_{n+1}^2\sin^2\psi=\big(x_1^2+\dots+x_n^2\big)\cos^2\psi
\;
\text{and}
\;
x_{n+1}>0\big\}
\\
\phantom{C}
= \{(x_1r\sin\psi,\ldots,x_nr\sin\psi,r\cos\psi):(x_1,\ldots,x_n)\in S^{n-1}
\;
\text{and}
\;
r>0\}.
\end{gather*}
Therefore, any parametrization of the unit sphere $S^{n-1}$ will give a~natural parametrization of the cone~$C$.
We choose the stereographic projection coordinates from the north pole:
\begin{gather*}
(x_1,\ldots,x_n)=\left(\frac{2v_1}{\sum\limits_i v_i^2+1},\ldots,\frac{2v_{n-1}}{\sum\limits_i
v_i^2+1},\frac{\sum\limits_i v_i^2-1}{\sum\limits_i
v_i^2+1}\right)=\left(\frac{2\mathbf{v}}{|\mathbf{v}|^2+1},\frac{|\mathbf{v}|^2-1}{|\mathbf{v}|^2+1}\right),
\end{gather*}
where $\mathbf{v}=(v_1,\ldots,v_{n-1})$.
This def\/ines a~coordinate system $\phi:\mathbb{R}^{n-1}\times\mathbb{R}^+\to C$ for the cone~$C$~by
\begin{gather}
\label{ConeParam}
\phi(\mathbf{v},r)=\left(\frac{2\mathbf{v}}{|\mathbf{v}|^2+1}r\sin\psi,\frac{|\mathbf{v}|^2-1}{|\mathbf{v}|^2+1}r\sin\psi,r\cos\psi\right),
\end{gather}
with inverse
\begin{gather}
\label{InverseConeParam}
\phi^{-1}(x_1,\ldots,x_{n+1}) = \left(\frac{x_1}{r\sin\psi-x_n},\ldots,\frac{x_{n-1}}{r\sin\psi-x_n},r\right),
\end{gather}
where $r:=\big(x_1^2+\dots+x_{n+1}^2\big)^{1/2}$.
Now, the metric~$g$ on~$C$ is the one inherited from the ambient Euclidean space $\mathbb{R}^{n+1}$.
In this coordinate system we get
\begin{gather*}
g=\frac{4r^2\sin^2\psi}{(|\mathbf{v}|^2+1)^2}\sum\limits_{i=1}^{n-1} \mathrm{d} v_i\otimes\mathrm{d} v_i+\mathrm{d}
r\otimes \mathrm{d} r.
\end{gather*}
The Christof\/fel symbols corresponding to this metric are
\begin{gather*}
\Gamma^k_{ij}=
\begin{cases}
2\dfrac{v_k\delta_{ij}-v_i\delta_{jk}-v_j\delta_{ki}}{|\mathbf{v}|^2+1} ,& i\neq n,\quad j\neq n,\quad k\neq n,
\\
-\dfrac{4r\sin^2\psi}{(|\mathbf{v}|^2+1)^2}\delta_{ij}+\dfrac{\delta_{ik}+\delta_{jk}}{r} ,&
\text{one and only one of~$i$,~$j$ or~$k$ is equal to~$n$},
\\
0,&
\text{else}.
\end{cases}
\end{gather*}
The geodesic equations are then
\begin{gather}
\label{ConeGeo1}
\ddot{\mathbf{v}}+\frac{2|\dot{\mathbf{v}}|^2\mathbf{v}-4(\mathbf{v}\cdot\dot{\mathbf{v}})\dot{\mathbf{v}}}{|\mathbf{v}|^2+1}
+\frac{2\dot{r}\dot{\mathbf{v}}}{r} = 0,
\\
\label{ConeGeo2}
\ddot{r}-\frac{4r\sin^2\psi|\dot{\mathbf{v}}|^2}{(|\mathbf{v}|^2+1)^2} = 0.
\end{gather}

The most important result of this section is the following.

\begin{theorem}
\label{TheoremConeGeo}
Let $\mathbf{r}:I\to\dot{\mathbb{R}}^{n+1}$ be a~non-colliding curve, where $n\geq 2$.
Then, $\mathbf{r}$ is a~geodesic on an~$n$-dimensional cone~$C$ if and only if $\mathbf{r}$ is a~geodesic on
a~$2$-dimensional cone $D\subseteq C$ of the same aperture.
\end{theorem}

\begin{proof}
We f\/irst show that if $\mathbf{r}:I\to\dot{\mathbb{R}}^{n+1}$ is a~non-colliding geodesic on~$C$, then $\mathbf{r}(I)$
is in a~2-dimensional cone $D\subseteq C$.
By def\/inition, we have to show that the curve $\alpha:I\to\mathbb{R}^{n+1}$ def\/ined~by
$\alpha(t):=\mathbf{r}(t)/|\mathbf{r}(t)|$ lies on a~f\/ixed 2-dimensional plane.
To do that, it suf\/f\/ices to show that $\{\dot{\alpha},\ddot{\alpha}\}$ is everywhere linearly independent while
$\{\dot{\alpha},\ddot{\alpha},\dddot{\alpha}\}$ is everywhere linearly dependent (see  \cite[Chapter~7, Part~B, Theorem~5]{Spivak}).

To show that $\{\dot{\alpha},\ddot{\alpha}\}$ is everywhere linearly independent, suppose that at some point $t_0\in I$
we have $\ddot{\alpha}=\lambda\dot{\alpha}$ for some $\lambda\in\mathbb{R}$.
Since $\alpha\cdot\alpha=1$, we get $\alpha\cdot\ddot{\alpha} = \lambda\alpha\cdot \dot{\alpha} = 0$.
Taking the second derivative on both sides of $\alpha\cdot\alpha=1$, we then f\/ind $\dot{\alpha}\cdot\dot{\alpha}=0$ at
$t_0$, and so $\dot{\mathbf{r}}(t_0)$ is parallel to $\mathbf{r}(t_0)$.
Hence, we can form a~colliding curve $\tilde{\mathbf{r}}(t):=\mathbf{r}(t_0)+(t-t_0)\dot{\mathbf{r}}(t_0)$, with
$\tilde{\mathbf{r}}(t_0)=\mathbf{r}(t_0)$ and $\dot{\tilde{\mathbf{r}}}(t_0)=\dot{\mathbf{r}}(t_0)$.
But $\tilde{\mathbf{r}}$ is solution to the geodesic equations~\eqref{ConeGeo1} and~\eqref{ConeGeo2}, so by uniqueness
we have $\mathbf{r}=\tilde{\mathbf{r}}$.
This contradicts the assumption that $\mathbf{r}$ is non-colliding.

Now, to show that $\{\dot{\alpha},\ddot{\alpha},\dddot{\alpha}\}$ is everywhere linearly dependent, we will show an
explicit non-trivial linear dependence.
In the parametrization $(\mathbf{v},r)$ of the cone~\eqref{ConeParam}, we have
\begin{gather*}
\alpha=\left(\frac{2\mathbf{v}}{|\mathbf{v}|^2+1}\sin\psi,\frac{|\mathbf{v}|^2-1}{|\mathbf{v}|^2+1}\sin\psi,\cos\psi\right).
\end{gather*}
To compute the derivatives $\dot{\alpha}$, $\ddot{\alpha}$, $\dddot{\alpha}$, we insert the geodesic
equations~\eqref{ConeGeo1} and~\eqref{ConeGeo2} to eliminate all second derivatives of $(\mathbf{v},r)$.
We obtain
\begin{gather*}
\dot{\alpha} =
\frac{2\sin\psi}{(|\mathbf{v}|^2+1)^2}
\big(\big(|\mathbf{v}|^2+1\big)\dot{\mathbf{v}}-2(\mathbf{v}\cdot\dot{\mathbf{v}})\mathbf{v},2\mathbf{v}\cdot\dot{\mathbf{v}},0\big),
\\
\ddot{\alpha} =
\frac{8\sin\psi}{(|\mathbf{v}|^2+1)^3}
\left(-|\dot{\mathbf{v}}|^2\mathbf{v}
+\frac{\dot{r}}{2r}\big(|\mathbf{v}|^2+1\big)
\big(2(\mathbf{v}\cdot\dot{\mathbf{v}})\mathbf{v}-\big(|\mathbf{v}|^2+1\big)\dot{\mathbf{v}}\big),\right.
\\
\left.\phantom{\ddot{\alpha}=}
{}-\frac{1}{2}|\dot{\mathbf{v}}|^2\big(|\mathbf{v}|^2-1\big)-\frac{\dot{r}}{r}\big(|\mathbf{v}|^2+1\big)\mathbf{v}\cdot\dot{\mathbf{v}},0\right),
\\
\dddot{\alpha}=
\frac{2\sin\psi}{(|\mathbf{v}|^2+1)^2}
\left(\frac{6\dot{r}^2}{r^2}-\frac{4(1+2\sin^2\psi)|\dot{\mathbf{v}}|^2}{(|\mathbf{v}|^2+1)^2}\right)
\big(\big(|\mathbf{v}|^2+1\big)\dot{\mathbf{v}}-2(\mathbf{v}\cdot\dot{\mathbf{v}})\mathbf{v},
2\mathbf{v}\cdot\dot{\mathbf{v}},0\big)
\\
\phantom{\dddot{\alpha}=}
{}-\frac{48\dot{r}\sin\psi}{r(|\mathbf{v}|^2+1)^3}\left(-|\dot{\mathbf{v}}|^2\mathbf{v},
-\frac{1}{2}|\dot{\mathbf{v}}|^2\big(|\mathbf{v}|^2-1\big),0\right).
\end{gather*}
We then f\/ind that these expressions satisfy the relation
\begin{gather*}
\left(\frac{4|\dot{\mathbf{v}}|^2(1+2\sin^2\psi)}{(|\mathbf{v}|^2+1)^2}+\frac{6\dot{r}^2}{r^2}\right)\dot{\alpha}
+\frac{6\dot{r}}{r}\ddot{\alpha}+\dddot{\alpha}=0.
\end{gather*}
So~$\alpha$ is contained in a~2-dimensional plane and hence $\mathbf{r}$ is contained in a~2-dimensional cone~$D$.
Moreover, $D\subseteq C$ for the following reason.
Let~$Q$ be the af\/f\/ine 2-dimensional plane gene\-ra\-ting~$D$, and let~$P$ be the af\/f\/ine hyperplane generating~$C$.
Since $\mathbf{r}$ is non-colliding, we can f\/ind three distinct points in $\alpha(I)\subseteq Q\cap S^n$.
But $Q\cap S^n$ is a~circle (Proposition~\ref{PlaneCapCircle}), so we have three non-collinear points of $Q\cap P$.
Since an af\/f\/ine 2-dimensional plane is uniquely specif\/ied by three non-collinear points, we have $Q\subseteq P$,
whence $D\subseteq C$.

We will now show that $\mathbf{r}$ is a~geodesic on~$D$ and that~$D$ has aperture~$\psi$.
Still assuming the parametrization~\eqref{ConeParam} for~$C$, we have
\begin{gather*}
P=\{(x_1,\ldots,x_n,\cos\psi):x_i\in\mathbb{R}\},
\end{gather*}
whence
\begin{gather*}
Q = (0,\ldots,0,\cos\psi)+ \mathbf{a}+\mathrm{span}\{\mathbf{w}_1,\mathbf{w}_2\},
\end{gather*}
for some $\mathbf{a},\mathbf{w}_1,\mathbf{w}_2\in\{(x_1,\ldots,x_n,0):x_i\in\mathbb{R}\}$.
But note that in the parametrization~\eqref{ConeParam},~$C$ has an $\mathrm{SO}(n)$ symmetry in its f\/irst~$n$
components.
Hence, we may assume that
\begin{gather*}
Q=\{(x_1,x_2,0,\ldots,0,a,\cos\psi):x_1,x_2\in\mathbb{R}\},
\end{gather*}
for some $a\geq 0$.
Since $Q\cap S^n$ contains more than one point, we have $a^2<1-\cos^2\psi=\sin^2\psi$, or equivalently,
$a=\cos\varphi\sin\psi$ for some $\varphi\in(0,\pi/2]$.
Therefore,
\begin{gather*}
D = \big\{r(x_1,x_2,0,\ldots,0,\cos\varphi\sin\psi,\cos\psi):r>0,x_1^2+x_2^2+\cos^2\varphi\sin^2\psi+\cos^2\psi=1\big\}
\\
\phantom{D}
 = \big\{r(x_1,x_2,0,\ldots,0,\cos\varphi\sin\psi,\cos\psi):r>0,x_1^2+x_2^2=\sin^2\varphi\sin^2\psi\big\}
\\
\phantom{D}
 = \big\{(r\cos\theta\sin\varphi\sin\psi,r\sin\theta\sin\varphi\sin\psi,0,\ldots,0,r\cos\varphi\sin\psi,r\cos\psi):r>0,\theta\in\mathbb{R}\big\},
\end{gather*}
and so
\begin{gather*}
\mathbf{r}(t)=(r(t)\cos\theta(t)\sin\varphi\sin\psi,r(t)\sin\theta(t)\sin\varphi\sin\psi,0,\ldots,0,r(t)\cos\varphi\sin\psi,r(t)\cos\psi)
\end{gather*}
for some functions $\theta:I\to\mathbb{R}$ and $r:I\to\mathbb{R}^+$.
Now, using the inverse transformation~\eqref{InverseConeParam}, we can express $\mathbf{r}$ in the coordinates
$(\mathbf{v},r)=(v_1,\ldots,v_{n-1},r)$ parametrizing~$C$.
We f\/ind
\begin{gather*}
(\mathbf{v},r)=\left(\cos\theta\cot\frac{\varphi}{2},\sin\theta\cot\frac{\varphi}{2},0,\ldots,0,r\right).
\end{gather*}
Assuming that $\mathbf{r}$ has this form, the geodesic equations~\eqref{ConeGeo1} and~\eqref{ConeGeo2} for $\mathbf{r}$
on~$C$ are equivalent~to
\begin{gather*}
\cos\varphi=0,
\\
\ddot{\theta}+2\dot{r}\dot{\theta}/r=0,
\\
\ddot{r}-r\dot{\theta}^2\sin^2\psi=0.
\end{gather*}
The f\/irst equation shows that~$D$ is the 2-dimensional cone of aperture~$\psi$ given~by
\begin{gather}
\label{2DCone}
D=\{(r\cos\theta\sin\psi,r\sin\theta\sin\psi,0,\ldots,0,r\cos\psi):r>0,\theta\in\mathbb{R}\},
\end{gather}
and the last two equations are precisely the geodesic equations for $\mathbf{r}$ on~$D$.

Conversely, this also shows that any geodesic of~$D$ is a~geodesic of~$C$.
Hence, by the $\mathrm{SO}(n)$ symmetry of~$C$ in its f\/irst~$n$ components (in~\eqref{ConeParam}) we get that any
geodesic on a~2-dimensional cone of aperture~$\psi$ embedded in~$C$ is a~geodesic on~$C$.
\end{proof}

The following two corollaries will be important for the next section, when we will investigate the motion of a~charged
particle in the f\/ield of Yang's monopole.

\begin{corollary}
\label{GeoSimpli}
Let $\mathbf{r}:I\to\dot{\mathbb{R}}^{n+1}$ be a~non-colliding curve, where $n\geq 2$.
If $\mathbf{r}$ lies on a~$2$-dimensional cone~$D$ of aperture~$\psi$ and $\ddot{\mathbf{r}}$ is orthogonal to
$\mathbf{r}$ and $\dot{\mathbf{r}}$, then $\mathbf{r}$ is a~geodesic on~$D$ and hence on the unique~$n$-dimensional cone
of aperture~$\psi$ containing~$D$.
\end{corollary}
\begin{proof}
By rotation symmetry, we may assume that the $2$-dimensional cone~$D$ on which $\mathbf{r}$ lies is of the
form~\eqref{2DCone}, so that
\begin{gather*}
\mathbf{r}(t)=(r(t)\cos\theta(t)\sin\psi,r(t)\sin\theta(t)\sin\psi,0,\ldots,0,r(t)\cos\psi),
\end{gather*}
for some functions $r:I\to\mathbb{R}^+$ and $\theta:I\to\mathbb{R}$.
This def\/ines a~space curve $\mathbf{p}:I\to\mathbb{R}^3$ by taking the three non-zero components of $\mathbf{r}$.
We have that $\mathbf{p}$ lies on the cone $\tilde{D}\subseteq\mathbb{R}^3$ of aperture~$\psi$ directed along $(0,0,1)$,
and $\ddot{\mathbf{p}}$ is orthogonal to $\mathbf{p}$ and $\dot{\mathbf{p}}$.
Hence, $\ddot{\mathbf{p}}$ is always normal to the surface of the cone $\tilde{D}$, so $\mathbf{p}$ is a~geodesic on
$\tilde{D}$.
Therefore, $\mathbf{r}$ is a~geodesic on~$D$.
Now, Proposition~\ref{UniqueConeEmbeding} shows that there exists a~unique~$n$-dimensional cone~$C$ of aperture~$\psi$
containing~$D$, and then Theorem~\ref{TheoremConeGeo} shows that $\mathbf{r}$ is a~geodesic on that cone.
\end{proof}

\begin{corollary}
\label{UniqueConeGeo}
Let $\mathbf{r}:I\to\dot{\mathbb{R}}^n$ be a~non-colliding geodesic on a~$k$-dimensional cone~$C$.
Then, $C$~is the unique~$k$-dimensional cone containing $\mathbf{r}(I)$.
\end{corollary}
\begin{proof}
Let~$D$ be any other~$k$-dimensional cone on which $\mathbf{r}$ is a~geodesic.
We want to show that $D=C$.
First suppose the case $k=2$ has been proved.
By Theorem~\ref{TheoremConeGeo}, $\mathbf{r}$ is a~non-colliding geodesic on a~2-dimensional cone $\tilde{C}\subseteq C$
of the same aperture as~$C$, and also on a~2-dimensional cone $\tilde{D}\subseteq D$ of the same aperture as~$D$.
Hence, $\tilde{C}=\tilde{D}$, and this cone has the same aperture~$\psi$ as~$C$ and~$D$.
Proposition~\ref{UniqueConeEmbeding} shows that~$C$ is the unique cone of aperture~$\psi$ containing
$\tilde{C}=\tilde{D}$, and the same is true for~$D$, so we have $C=D$.
We may thus assume that~$C$ and~$D$ are $2$-dimensional.

\looseness=1
Now, since $\mathbf{r}$ is non-colliding, we can f\/ind three points in $\mathbf{r}(I)\subseteq C$ such that no two of
them are collinear with the origin.
Using Proposition~\ref{PlaneCapCircle}, we f\/ind that the radial projection of these points on the unit sphere gives 3
distinct points on a~circle.
Hence, we get 3 non-collinear points on the 2-dimensional af\/f\/ine plane~$P$ generating~$C$.
Now, the 2-dimensional af\/f\/ine plane~$Q$ generating~$D$ must also contain these 3 non-collinear points.
Since an af\/f\/ine 2-dimensional plane in $\mathbb{R}^n$ is uniquely def\/ined by 3 non-collinear points, we have $Q=P$, and
hence $D=C$.
\end{proof}

\section{Particle motion in Yang's monopole and geodesics on cones}
\label{Section8}

In this section we investigate the solutions to the Poincar\'e problem in $\dot{\mathbb{R}}^5$.
That is, we solve the equations~\eqref{SU2EOM1},~\eqref{SU2EOM2} and~\eqref{SU2EOM3} for the motion of a~charged
particle in the f\/ield Yang's $\mathrm{SU}(2)$ monopole at the origin of Euclidean space $\mathbb{R}^5$.

Let us denote a~solution to the Poincar\'e problem in $\dot{\mathbb{R}}^5$ by a~pair $(\mathbf{r},\mathbf{e})$ of curves
$\mathbf{r}:I\to\dot{\mathbb{R}}^5$ and $\mathbf{e}:I\to\mathbb{R}^3$, for some open interval~$I$.
More precisely, $(\mathbf{r},\mathbf{e})$ is a~solution if the curve $(\mathbf{u},r):I\to\mathbb{R}^4\times\mathbb{R}^+$
obtained by expressing $\mathbf{r}$ in the stereographic projection coordinates from the north pole~\eqref{StereoNorth}
together with the curve $\mathbf{e}=(e^1,e^2,e^2):I\to\mathbb{R}^3$ satisfy the equations of
motion~\eqref{SU2EOM1},~\eqref{SU2EOM2} and~\eqref{SU2EOM3} for all $t\in I$.

Our main goal is to show that for every solution $(\mathbf{r},\mathbf{e})$ there is a~4-dimensional cone with vertex
at the origin of $\mathbb{R}^5$ on which $\mathbf{r}$ is a~geodesic.
Note that this fact together with Theorem~\ref{TheoremConeGeo} show that $\mathbf{r}$ is also a~geodesic on
a~2-dimensional cone, as was the case for every solution to the Poincar\'e problem in $\dot{\mathbb{R}}^3$.
It is quite remarkable that although we are dealing with a~non-Abelian monopole and hence far more intricate equations
of motion, the space of solution is almost identical to the one describing motion of a~particle in the simpler Abelian
Dirac monopole.

For Dirac's monopole, the hard part of the proof is to f\/ind an explicit expression for the direction
$\mathbf{L}\in\dot{\mathbb{R}}^3$ of the cone.
Once we have $\mathbf{L}$, it is very easy to see that $\mathbf{r}$ is at a~constant angle from $\mathbf{L}$ and that
$\ddot{\mathbf{r}}$ is always normal to the surface of the cone.
For Yang's monopole, we will also f\/ind a~vector $\mathbf{L}\in\dot{\mathbb{R}}^5$ for which $\mathbf{r}$ is at
a~constant angle, which will then imply that $\mathbf{r}$ lies on a~4-dimensional cone~$C$.
However, the proof that $\mathbf{r}$ is a~geodesic on~$C$ is more tricky.
We will show that $\ddot{\mathbf{r}}$ is orthogonal to $\dot{\mathbf{r}}$ and $\mathbf{r}$, as we did for Dirac's
monopole, but in $\mathbb{R}^5$ this fact is not suf\/f\/icient to infer that $\ddot{\mathbf{r}}$ is normal to the surface
of the cone.
To complete the proof we will need to use some non-trivial conclusions of the preceding section on higher dimensional
cones, namely, Corollaries~\ref{GeoSimpli} and~\ref{UniqueConeGeo}.

Now, looking at the equations of motion~\eqref{SU2EOM1},~\eqref{SU2EOM2} and~\eqref{SU2EOM3} for the Poincar\'e problem
in~$\dot{\mathbb{R}}^5$, we immediately see that a colliding curve is a solution if and only if it has constant speed, as was the case for Dirac's monopole. Since a~constant-speed colliding curve is a~geodesic of many cones, we may exclude these trivial solutions from our discussion.
The main result of our paper is the following.

\begin{theorem}
\label{TheoremYangMonopole}
Let $(\mathbf{r},\mathbf{e})$ be a~solution to the Poincar\'e problem in $\dot{\mathbb{R}}^5$.
If $\mathbf{r}$ is non-colliding, then $\mathbf{r}$ is a~geodesic on the $4$-dimensional cone directed along the constant
vector
\begin{gather}
\label{ConstantVector}
\mathbf{L}:=
\left(\frac{(|\mathbf{e}|^2-4r^2(\mathbf{A}\dot{\mathbf{u}}\cdot\mathbf{e}))\mathbf{u}+2r^2\mathbf{E}\dot{\mathbf{u}}}{2(|\mathbf{u}|^2+1)},
\frac{2r^2(\mathbf{A}\dot{\mathbf{u}}\cdot\mathbf{e})}{|\mathbf{u}|^2+1}
+\frac{|\mathbf{e}|^2}{4}\frac{|\mathbf{u}|^2-1}{|\mathbf{u}|^2+1}\right)\in\mathbb{R}^5
\end{gather}
and of aperture~$\psi$ given~by
\begin{gather*}
\cos\psi = \frac{|\mathbf{e}|}{2}\left(\frac{|\mathbf{e}|^2}{4}+ \frac{4
r^4|\dot{\mathbf{u}}|^2}{(|\mathbf{u}|^2+1)^2}\right)^{-1/2}.
\end{gather*}
\end{theorem}

\begin{proof}
We f\/irst show that $\mathbf{r}$ is a~geodesic on \emph{some} 4-dimensional cone.
By Corollary~\ref{GeoSimpli}, it suf\/f\/ices to show that $\mathbf{r}$ lies on a~2-dimensional cone and $\ddot{\mathbf{r}}$
is orthogonal to $\mathbf{r}$ and $\dot{\mathbf{r}}$.
In the stereographic projection coordinates from the north pole $(\mathbf{u},r)=(u_1,\ldots,u_4,r)$, we have
\begin{gather*}
\mathbf{r}=\left(\frac{2r\mathbf{u}}{|\mathbf{u}|^2+1},r\frac{|\mathbf{u}|^2-1}{|\mathbf{u}|^2+1}\right),
\end{gather*}
and a~straightforward computation shows that
\begin{gather*}
\ddot{\mathbf{r}}\cdot\mathbf{r} = \ddot{r}r-\frac{4r^2|\dot{\mathbf{u}}|^2}{(|\mathbf{u}|^2+1)^2},
\qquad
\ddot{\mathbf{r}}\cdot\dot{\mathbf{r}} =
\ddot{r}\dot{r}+4\frac{\dot{r}r|\dot{\mathbf{u}}|^2
+r^2\dot{\mathbf{u}}\cdot\ddot{\mathbf{u}}}{(|\mathbf{u}|^2+1)^2}
-8\frac{r^2|\dot{\mathbf{u}}|^2(\mathbf{u}\cdot\dot{\mathbf{u}})}{(|\mathbf{u}|^2+1)^3}.
\end{gather*}
By inserting the equations of motion~\eqref{SU2EOM1} and~\eqref{SU2EOM2}, we immediately get that these two expressions
are equal to zero.

Now, to show that $\mathbf{r}$ lies on a~2-dimensional cone, we will follow an approach very similar to the one in the
proof of Theorem~\ref{TheoremConeGeo}.
That is, let $\alpha:I\to\mathbb{R}^{5}$ be def\/ined by $\alpha(t):=\mathbf{r}(t)/|\mathbf{r}(t)|$.
It suf\/f\/ices to show that $\alpha(I)$ is contained in a~2-dimensional af\/f\/ine plane, or equivalently, that
$\{\dot{\alpha},\ddot{\alpha}\}$ is everywhere linearly independent while
$\{\dot{\alpha},\ddot{\alpha},\dddot{\alpha}\}$ is everywhere linearly dependent.
The proof that $\{\dot{\alpha},\ddot{\alpha}\}$ is everywhere linearly independent is exactly the same as the one in the
proof of Theorem~\ref{TheoremConeGeo}.
Now, we will show an explicit non-trivial linear dependence of $\{\dot{\alpha},\ddot{\alpha},\dddot{\alpha}\}$.
First, we have
\begin{gather*}
\alpha = \left(\frac{2\mathbf{u}}{|\mathbf{u}|^2+1},\frac{|\mathbf{u}|^2-1}{|\mathbf{u}|^2+1}\right),
\qquad
\dot{\alpha} =
\left(\frac{2\dot{\mathbf{u}}}{|\mathbf{u}|^2+1}-\frac{4(\mathbf{u}\cdot\dot{\mathbf{u}})\mathbf{u}}{(|\mathbf{u}|^2+1)^2},
\frac{4\mathbf{u}\cdot\dot{\mathbf{u}}}{(|\mathbf{u}|^2+1)^2}\right).
\end{gather*}
We then compute the second and third derivatives of~$\alpha$ by inserting the equations of
motion~\eqref{SU2EOM1},~\eqref{SU2EOM2} and~\eqref{SU2EOM3} to eliminate all second and higher derivatives of
$(\mathbf{u},r)$.
We f\/ind
\begin{gather*}
\ddot{\alpha} =
\left(\frac{\mathbf{E}\dot{\mathbf{u}}-2(\mathbf{A}\dot{\mathbf{u}}\cdot\mathbf{e})\mathbf{u}}{r^2(|\mathbf{u}|^2+1)}
-\frac{8|\dot{\mathbf{u}}|^2\mathbf{u}}{(|\mathbf{u}|^2+1)^3}
+\frac{4\dot{r}}{r}\left(\frac{2(\mathbf{u}\cdot\dot{\mathbf{u}})\mathbf{u}}{(|\mathbf{u}|^2+1)^2}-\frac{\dot{\mathbf{u}}}{|\mathbf{u}|^2+1}\right),
\right.
\\
\phantom{\ddot{\alpha}=}
\left.
\frac{2\mathbf{A}\dot{\mathbf{u}}\cdot\mathbf{e}}{r^2(|\mathbf{u}|^2+1)}-\frac{4|\dot{\mathbf{u}}|^2(|\mathbf{u}|^2-1)}{(|\mathbf{u}|^2+1)^3}
-\frac{8\dot{r}\mathbf{u}\cdot\dot{\mathbf{u}}}{r(|\mathbf{u}|^2+1)^2}\right),
\\
\dddot{\alpha} = \left(\frac{12
r\dot{r}(2(\mathbf{A}\dot{\mathbf{u}}\cdot\mathbf{e})\mathbf{u}-\mathbf{E}\dot{\mathbf{u}})-|\mathbf{e}|^2\dot{\mathbf{u}}}{2r^4(|\mathbf{u}|^2+1)}
+\frac{|\mathbf{e}|^2(\mathbf{u}\cdot\dot{\mathbf{u}})\mathbf{u}}{r^4(|\mathbf{u}|^2+1)^2}
+\frac{24|\dot{\mathbf{u}}|^2(2\dot{r}\mathbf{u}-r\dot{\mathbf{u}})}{r(|\mathbf{u}|^2+1)^3}\right.
\\
\phantom{\dddot{\alpha}=}
{}+\frac{48|\dot{\mathbf{u}}|^2(\mathbf{u}\cdot\dot{\mathbf{u}})\mathbf{u}}{(|\mathbf{u}|^2+1)^4}
+\frac{12\dot{r}^2}{r^2}\left(\frac{\dot{\mathbf{u}}}{|\mathbf{u}|^2+1}
-\frac{2(\mathbf{u}\cdot\dot{\mathbf{u}})\mathbf{u}}{(|\mathbf{u}|^2+1)^2}\right),
-\frac{12\dot{r}\mathbf{A}\dot{\mathbf{u}}\cdot\mathbf{e}}{r^3(|\mathbf{u}|^2+1)}
\\
\phantom{\dddot{\alpha}=}
\left.-\frac{|\mathbf{e}|^2\mathbf{u}\cdot\dot{\mathbf{u}}}{r^4(|\mathbf{u}|^2+1)^2}
+\frac{24\dot{r}|\dot{\mathbf{u}}|^2(|\mathbf{u}|^2-1)}{r(|\mathbf{u}|^2+1)^3}
-\frac{48|\dot{\mathbf{u}}|^2\mathbf{u}\cdot\dot{\mathbf{u}}}{(|\mathbf{u}|^2+1)^4}
+\frac{24\dot{r}^2\mathbf{u}\cdot\dot{\mathbf{u}}}{r^2(|\mathbf{u}|^2+1)^2}\right).
\end{gather*}
To perform these computations, it is useful to f\/irst derive the following identities
\begin{gather*}
\mathbf{u}\cdot\mathbf{E}\dot{\mathbf{u}} = (|\mathbf{u}|^2+1)\mathbf{A}\dot{\mathbf{u}}\cdot\mathbf{e},
\\
\mathbf{e}\mathbf{A}\mathbf{E} =-\frac{|\mathbf{e}|^2\mathbf{u}}{|\mathbf{u}|^2+1},
\\
\dot{\mathbf{E}}\dot{\mathbf{u}} =
2\frac{|\dot{\mathbf{u}}|^2\mathbf{E}\mathbf{u}+(\mathbf{u}\cdot\mathbf{E}\dot{\mathbf{u}})\dot{\mathbf{u}}
-(\mathbf{u}\cdot\dot{\mathbf{u}})\mathbf{E}\dot{\mathbf{u}}}{|\mathbf{u}|^2+1}.
\end{gather*}
We then get the simple expression
\begin{gather*}
\left(\frac{|\mathbf{e}|^2}{4r^2}+\frac{6\dot{r}^2}{r^2}+\frac{12|\dot{\mathbf{u}}|^2}{(|\mathbf{u}|^2+1)^2}\right)\dot{\alpha}
+\frac{6\dot{r}}{r}\ddot{\alpha}+\dddot{\alpha}=0.
\end{gather*}
Hence, $\{\dot{\alpha},\ddot{\alpha},\dddot{\alpha}\}$ is everywhere linearly dependent, so $\mathbf{r}$ lies on
a~2-dimensional cone, and by Corollary~\ref{GeoSimpli}, $\mathbf{r}$ is a~geodesic on a~$4$-dimensional cone~$C$.

We will now give an explicit expression for the cone~$C$.
By dif\/ferentiating the vector $\mathbf{L}$ given by~\eqref{ConstantVector} and inserting the equations of motion, we see
that $\mathbf{L}$ is constant.
Moreover,
\begin{gather*}
\frac{\mathbf{r}\cdot\mathbf{L}}{|\mathbf{r}||\mathbf{L}|}=\frac{|\mathbf{e}|^2}{4|\mathbf{L}|} =
\frac{|\mathbf{e}|}{2}\left(\frac{|\mathbf{e}|^2}{4}+ \frac{4
r^4|\dot{\mathbf{u}}|^2}{(|\mathbf{u}|^2+1)^2}\right)^{-1/2},
\end{gather*}
which is also constant.
Since $\mathbf{r}$ is non-colliding we have $\dot{\mathbf{u}}\neq 0$, so this expression can be written as the cosine of
some angle $\psi\in(0,\pi/2]$.
Hence, $\mathbf{r}$ lies on the 4-dimensional cone of aperture~$\psi$ directed along $\mathbf{L}$.
But we showed that $\mathbf{r}$ is a~geodesic on some 4-dimensional cone~$C$, and by Corollary~\ref{UniqueConeGeo},~$C$
is the unique 4-dimensional cone containing $\mathbf{r}(I)$.
Hence,~$C$ is the cone of aperture~$\psi$ directed along $\mathbf{L}$.
\end{proof}

By this theorem, the problem reduces to the geodesic equations on a~4-dimensional cone, which in turn reduces to the
geodesic equations on a~2-dimensional cone by Theorem~\ref{TheoremConeGeo}.
Geodesics on 2-dimensional cones were discussed in Section~\ref{Section2}.

\subsection*{Acknowledgements}

The author is grateful to Professor Niky Kamran for his constant guidance and invaluable suggestions.
The author would also like to thank the anonymous referees who provided helpful comments, corrections and reference
suggestions.
This work was supported by the NSERC USRA program, grant number RGPIN 105490-2011.

\pdfbookmark[1]{References}{ref}
\LastPageEnding

\end{document}